\documentclass[11pt,reqno]{amsart}
\usepackage[centertags]{amsmath}
\usepackage{amsfonts}
\usepackage{amssymb}
\usepackage{amsthm}
\usepackage{mathrsfs}
\newtheorem{theorem}{Theorem}
\newtheorem{corollary}{Corollary}
\newtheorem{lemma}{Lemma}
\theoremstyle{definition}
\newtheorem{remark}{Remark}
\begin{document}
%
%
\title{Cooperative Stochastic Control for Optical Beam Tracking}
\thanks{This research was supported by the Army Research Office
under ODDR\&E MURI01 Program Grant No. DAAD19-01-1-0465 to the
Center for Communicating Networked Control Systems (through Boston
University), and by the National Science Foundation under Grant No.
ECS 0636613.}
%
%
\author{Arash Komaee}%
\dedicatory{\rm Department of Electrical and Computer Engineering\\
University of Maryland, College Park, MD 20742, USA\\}
\thanks{E-mail address: \tt akomaee@eng.umd.edu}
%
%
\begin{abstract}
Maintaining optical alignment between stations of a free-space
optical link requires a persistent beam tracking operation. This is
achieved using a position-sensitive photodetector at each station
which measures the azimuth and elevation of tracking error. A
pointing assembly adjusts the heading of transceivers according to
measurement of the tracking error. The measurement at each receiver
also depends on the pointing error of the opposite transmitter,
therefore a cooperative beam tracking system can be regarded as two
dynamically coupled subsystems. We developed a stochastic model for
a cooperative beam tracking system to get insight into solution of
an associated optimal control problem with goal of maximizing the
flow of optical power between the stations.
\end{abstract}
%
%
\maketitle
%
%
\section{Introduction}
In free-space optical communication using narrow laser beams, it is
required to maintain the alignment of transmitter and receiver
stations in spite of their relative motion. This relative motion
might be caused by the mobile nature of the stations, mechanical
vibration, or accidental shocks. Prior to data transmission begins,
coarse alignment is achieved through two operations:
\textit{pointing} and \textit{spatial acquisition}. Pointing is the
act of aiming the transmitted beam toward the receiver within an
acceptable accuracy. The purpose of spatial acquisition is to detect
the transmitter's beam and align the normal vector to the receiving
optical device with the direction of the impinging optical field.
Following a coarse alignment accomplished in the acquisition phase,
data transmission is established and simultaneously the operation of
\textit{cooperative beam tracking} is performed. This fine alignment
operation is intended to precisely compensate for persistent
relative motion of the stations. For detailed description of these
three operations we refer the reader
to~\cite{BOOK.Gagliardi.95,BOOK.Katzman.87.ch6}. In this study we
are concerned about cooperative beam tracking.

In a cooperative beam tracking system, the stations continuously
track the arrival direction of the incident beam and transmit their
beam back in that direction. A position-sensitive photodetector
(e.g. quadrant detector) in association with a focusing lens is
employed at each station to measure the azimuth and elevation
components of the error vector. The error vector is the displacement
of the beam's arrival direction with respect to the direction normal
to the receiving aperture. A servo-driven pointing assembly adjusts
the heading of the transmitting optical device (in azimuth and
elevation directions) according to the measured error vector.
Normally, a single pointing assembly is used to control the
direction of both receiving and transmitting optics, which are
installed on the same platform.

Cooperative tracking systems were already studied by other
researchers based on deterministic
models~\cite{BOOK.Gagliardi.95,PROC.Wei.88,ART.Marola.89}. These
models describe the position-sensitive photodetector by a
deterministic input-output relationship and employ deterministic
functions to characterize the relative motion of the stations. Using
such a deterministic model, under a proportional control law, Wei
and Gagliardi~\cite{PROC.Wei.88} evaluated the steady-state
performance of the system and Marola~\textit{et~al.} studied
stability properties~\cite{ART.Marola.89}.

For most applications, a stochastic model may be more appropriate
for relative motion of the stations. Moreover, the optical sensing
devices are usually described by stochastic
models~\cite{BOOK.Gagliardi.95,BOOK.Goodman.85}. This leads us to
suggest that a stochastic approach is more likely to lead to
rigorous analysis and design of cooperative tracking systems. In
this paper, we develop a stochastic model for cooperative tracking
systems and employ the model to analyze the system and study an
associated optimal control problem. Maximizing the flow of the
optical power between the stations shall be our criterion for
optimality.
%
%
\section{System Description}
We consider an optical transceiver comprised of a lens, a
position-sensitive photodetector, and a narrow laser source, all
installed on a rigid platform. The photodetector surface is
perpendicular to the lens axis and its center is placed at the focus
of the lens. The axes of the lens and the laser source are parallel
to the transceiver axis. The azimuth and elevation of the
transceiver axis can be controlled by means of a servo-driven
pointing assembly. A two-way optical link employs two transceivers
of this type in such a manner that each transceiver transmits its
optical beam toward the opposite station and receives the optical
beam from the opposite side. The optical beams are used for two
purposes: as a carrier of information and as a beacon assisting the
opposite station in its tracking operation. We assume that the
stations are subjected to relative motion.

In what follows, we distinguish the stations by superscripts~$a$
and~$b$ or~$i=a,b$ when referring to both stations. The dependence
on time will be shown by subscript~$t$. The two-dimensional
vector~$\theta_t^i$ denotes the azimuth and elevation angles of the
transceiver axis~$i$ with respect to some fixed coordinate system.
Similarly,~$\phi_t$ denotes the azimuth and elevation angles of the
line-of-sight (LOS) of the stations with respect to same coordinate
system. We define the tracking error of station~$i$ as
$\psi_t^i=\theta_t^i-\phi_t$ which is equivalent to the pointing
error for the opposite station.

We shall assume that the transmitted optical fields impinge the
receivers along the line-of-sight of the stations, regardless of the
pointing error of the transmitters. This implies that the received
optical field at station~$i$ strikes the receiving lens along the
error vector~$\psi_t^i$ with respect to the axis of the lens. The
validity of this assumption is clear for a spherical optical field.
For a Gaussian beam~\cite{BOOK.Davis.00}, which is the model used in
this study, we can show that the angle between LOS and the beam
arrival direction at the receiver depends on the third power of the
pointing error at the transmitter. Because the pointing error is
maintained small by means of feedback control, its third power can
be neglected at least for a linear model, which justifies our
assumption.

In contrast to the beam arrival direction, the intensity of the
received optical field substantially depends on the transmitter's
pointing error. Assume that station~$b$ transmits a circular
symmetric Gaussian beam with divergence angle~$2\bar{\psi}$ toward
station~$a$. Then, due to the pointing error~$\psi_t^b$ at
station~$b$, the instantaneous optical power received by station~$a$
at time~$t$ is reduced by a factor of
$\exp\left(-2\|\psi_t^b\|^2/\bar{\psi}^2\right)$. This attenuation
factor is obtained based on two assumption: first, the radius of the
optical beam is much larger than the receiver aperture, and second,
the distance between the stations is short enough to allow ignoring
the propagation delay.

The image of the received optical field over the surface of the
photodetector is a spot of light with a bell-shaped intensity
profile centered at~$y_t^i=f_c\psi_t^i$, where~$f_c$ is the focal
length of the lens~\cite{BOOK.Gagliardi.95}.
Let~$\Omega\left(r\right)$ be the intensity pattern of the spot of
light, where~$r$ is the position vector of a point on the surface of
the photodetector. Denote by~$P_t^a>0$ the total optical power
received by station~$a$ at time~$t$ in the absence of pointing
error~$\psi_t^b$. Then the optical intensity over the surface of
photodetector~$a$ is given~by
\begin{equation}\label{CooperativeOpticalIntensityErrorVector}
I_t^a\left(r\right)=P_t^a\exp\textstyle\left(-2\|\psi_t^b\|^2/\bar{\psi}^2\right)\Omega\left(r-f_c\psi_t^a\right)
\end{equation}
A similar expression can be obtained for~$I_t^b\left(r\right)$ by
flipping~$a$ and~$b$
in~\eqref{CooperativeOpticalIntensityErrorVector}. Note that the
displacement of the spot of light depends only on the error vector
of the same station, while the total received power depends on the
error vector of the opposite station.

Since~$\psi_t^i$ depends linearly on the displacement of the spot of
light, it can be estimated from the output of the position-sensitive
photodetector. This estimate is provided to a controller which
applies proper control signals to the pointing assembly in order to
drive the error vector to zero. Note that $\psi_t^a=\psi_t^b=0$ is
the ``unattainable'' goal of a cooperative tracking system. Under a
realistic condition, the objective of the system is to
maintain~$\|\psi_t^a\|$ and~$\|\psi_t^b\|$ as small as possible.
Since the axis of the laser source is parallel to the transceiver's
axis, $\psi_t^i=0$ implies that the station~$i$ transmits its
optical beam along LOS, which leads to maximum flow of the optical
power between the stations.

From the above description, we find out that a cooperative tracking
system consists of two dynamical subsystems coupled via their
measurement. The subsystems operate cooperatively in the sense that
a small pointing error at one station assists the other station by
increasing its received optical power which leads to more accurate
estimation of its error vector.
%
%
\section{The Model}
The model we use in this study is a two-station extension of the
single-station model in~\cite{PROC.Snyder.78}. We refer the reader
to that article for detailed description and justification of the
model. Without loss of generality, we assume the stations are
identical.

The pointing assembly is an electro-mechanical system with the input
vector~$u_t^i\in\mathbb{R}^2$ and the output
vector~$\theta_t^i\in\mathbb{R}^2$. The input and output vectors are
two-dimensional corresponding to the azimuth and elevation angles.
We model this system by the linear stochastic differential equation
\begin{equation}\label{CooperativePointingAssembly}
\begin{split}
dx_t^{p,i}&=A_t^px_t^{p,i}dt+B_t^pu_t^idt+D_t^pdw_t^{p,i} \\
\theta_t^i&=C_t^px_t^{p,i}
\end{split}
\end{equation}
where $x_t^{p,i}\in\mathbb{R}^{n_p}$ is the state vector,
$\{w_t^{p,i},t\geqslant0\}$ is a $m_p$-dimensional standard Wiener
process, and~$A_t^p$, $B_t^p$, $D_t^p$, and $C_t^p$ are uniformly
bounded matrices with proper dimensions. We assume that
$\{w_t^{p,a},t\geqslant0\}$ and $\{w_t^{p,b},t\geqslant0\}$ are
mutually independent.

Using a linear model for the pointing assembly is justified by the
fact that the system operates over small angles during the fine
control regime. In applications like intersatellite communication,
the relative motion consists of a large, deterministic component and
a small, stochastic term. Accordingly, the control law consists of a
deterministic, open-loop, coarse control and a small, closed-loop,
fine control. In this case, the nonlinear state equations describing
the system is linearized around the deterministic nominal
trajectory, which results in the time-varying
model~\eqref{CooperativePointingAssembly} for the fine control
regime.

We model~$\phi_t$ by a Gauss-Markov stochastic process described by
the state-space equations
\begin{equation}\label{CooperativeLineOfSight}
\begin{split}
dx_t^d&=A_t^dx_t^ddt+D_t^ddw_t^d \\
\phi_t&=C_t^dx_t^d
\end{split}
\end{equation}
with state vector~$x_t^d\in\mathbb{R}^{n_d}$, $m_d$-dimensional
standard Wiener process~$\{w_t^d,t\geqslant0\}$, and uniformly
bounded matrices $A_t^d$, $D_t^d$, and $C_t^d$ with proper
dimensions.

The position vector~$y_t^i=f_c\psi_t^i$ is a linear function
of~$x_t^{p,i}$ and~$x_t^d$, so we can combine
equations~\eqref{CooperativePointingAssembly}
and~\eqref{CooperativeLineOfSight} and write them in the compact
form
\begin{equation}\label{CooperativeStateSpaceModel}
\begin{split}
dx_t^i&=A_tx_t^idt+B_tu_t^idt+D_tdw_t^i \\
y_t^i&=C_tx_t^i
\end{split}
\end{equation}
with state vector~$x_t^i\in\mathbb{R}^n$ and $m$-dimensional
standard Wiener process~$\{w_t^i,t\geqslant0\}$, where~$n=n_p+n_d$
and~$m=m_p+m_d$. The initial state~$x_0^i$ is assumed to be a
Gaussian vector with mean~$\bar{x}_0^i$ and covariance
matrix~$\bar{\Sigma}_0^i$ and independent of~$\{w_t^a,t\geqslant0\}$
and~$\{w_t^b,t\geqslant0\}$.

We approximate the bell-shaped intensity
profile~$\Omega\left(r\right)$
in~\eqref{CooperativeOpticalIntensityErrorVector} by a Gaussian
function. Then, in terms of the state vectors~$x_t^a$ and~$x_t^b$,
the optical intensity~$I_t^a\left(r\right)$ can be expressed~as
\begin{equation}\label{CooperativeOpticalIntensityStateVector}
I_t^a\left(r\right)=P_t^a\exp\left(-\rho\|C_tx_t^b\|^2\right)\gamma_t\left(r,x_t^a\right)
\end{equation}
where~$\rho=2/(\bar{\psi}f_c)^2$ and~$\gamma_t\left(r,x\right)$ is
defined~as
\begin{equation}\label{CooperativeGaussianIntensityProfile}
\gamma_t\left(r,x\right)=\left(2\pi\right)^{-1}\left(\det{R_t}\right)^{-1/2}
\exp\left\{-\frac{1}{2}\left(r-C_tx\right)^T(R_t)^{-1}\left(r-C_tx\right)\right\}
\end{equation}
Here, $R_t=R_t^T$ is a $2\times2$ positive-definite matrix
describing the shape of the pattern. For a circular symmetric
pattern with constant radius $\varrho>0$ we have
$R_t={\varrho}I_{2\times2}$. A similar expression can be derived
for~$I_t^b\left(r\right)$ by exchanging~$a$ and~$b$. We remind that
in~\eqref{CooperativeOpticalIntensityStateVector} the propagation
delay is neglected.

We allow~$P_t^a$ and~$P_t^b$ to be nonnegative stochastic processes
with piecewise continuous sample paths and nonzero expectations to
model the random optical fade caused by atmospheric turbulence and
aerosols and the information-bearing signals modulating the optical
beams. Further, we assume that $\{P_t^a,t\geqslant0\}$ and
$\{P_t^b,t\geqslant0\}$ are mutually independent and independent
of~$x_0^i$ and $\{w_t^i,t\geqslant0\},~i=a,b$.

The position-sensitive photodetector is a photoelectron converter
whose surface is partitioned into small regions. The output of each
region counts the number of converted electrons regardless of their
location on the region. The photoelectron conversion rate depends
linearly on the optical power absorbed by the region. Generally, a
photoelectron converter is modeled by a Poisson process with a rate
proportional to the impinging optical
power~\cite{BOOK.Gagliardi.95,BOOK.Goodman.85}. In the present case,
where the optical power is a stochastic field, the output of each
region shall be modeled by a conditionally (doubly stochastic)
Poisson process.

Many practical beam tracking systems employ a quadrant detector, a
photodetector with a four-region partition, as their optical sensing
device. However, the low spatial resolution of the quadrant detector
can be improved using a finer partition. For instance, the authors
of~\cite{PROC.LeeAlexanderJeganathan.00} describe a beam tracking
system utilizing a photodetector with $512\times512$ pixels. In this
study, following~\cite{PROC.Snyder.78}, we use an infinite
resolution model for the photodetector. This idealized model
provides a reasonable approximation for high resolution
photodetectors. We also make another ideal assumption that the
surface of the photodetector is infinitely
large~\cite{PROC.Snyder.78}. This assumption is justified when the
photodetector area is significantly larger than the size and the
displacement of the spot of light. We believe that the control law
obtained from this idealized model provides a useful point of
departure for practical designs, even for low resolution or finite
area photodetectors.

We use a doubly stochastic space-time Poisson process to describe
the output of an infinite resolution
photodetector~\cite{PROC.Snyder.78}. The rate of this process is
assumed to be proportional to~$I_t^i\left(r\right)$.
From~\eqref{CooperativeOpticalIntensityStateVector}, the stochastic
rate associated to photodetector~$i$ can be expressed~as
\begin{equation}\label{CooperativePoissonRate}
\lambda_t\left(r,x_t^i,\mu_t^i\right)=\mu_t^i\gamma_t\left(r,x_t^i\right)
\end{equation}
where~$\mu_t^a$ and~$\mu_t^b$ are defined~as
\begin{equation}\label{CooperativeMUDefinition}
\begin{split}
\mu_t^a&=\nu_t^a\exp\bigl(-\rho\|C_tx_t^b\|^2\bigr)\\
\mu_t^b&=\nu_t^b\exp\bigl(-\rho\|C_tx_t^a\|^2\bigr)
\end{split}
\end{equation}
Here,~$\nu_t^i$ is defined as~$\nu_t^i={\eta}P_t^i$, where $\eta>0$
is the sensitivity of the photodetectors and is assumed to be a
constant. Note that $\nu_t^i$ inherits the statistical properties
of~$P_t^i$.

The space-time Poisson process, defined
over~$[0,\infty)\times\mathbb{R}^2$, characterizes the occurrence of
discrete events (e.g. release of a single electron) with a temporal
component $t\in[0,\infty)$ and a spatial
component~$r\in\mathbb{R}^2$. Let $\mathcal{T}$ and $\mathcal{S}$ be
Borel sets in $[0,\infty)$ and~$\mathbb{R}^2$ respectively and
$N^i\left(\mathcal{T}\times\mathcal{S}\right)$ denote the number of
points occurring in $\mathcal{T}\times\mathcal{S}$. Define the
random variable
\begin{equation*}
\Lambda^i\left(\mathcal{T}\times\mathcal{S}\right)
=\int_{\mathcal{T}\times\mathcal{S}}\lambda_t\left(r,x_t^i,\mu_t^i\right)dtdr
\end{equation*}
Then $N^i\left(\mathcal{T}\times\mathcal{S}\right)$ is a doubly
stochastic Poisson random variable with conditional probability
distribution
\begin{equation*}
\Pr\left\{N^i\left(\mathcal{T}\times\mathcal{S}\right)=n|\Lambda^i\left(\mathcal{T}\times\mathcal{S}\right)\right\}
=\frac{e^{-\Lambda^i\left(\mathcal{T}\times\mathcal{S}\right)}
\left(\Lambda^i\left(\mathcal{T}\times\mathcal{S}\right)\right)^n}{n!}
\end{equation*}
Moreover, for disjoint $\mathcal{T}_1\times\mathcal{S}_1$ and
$\mathcal{T}_2\times\mathcal{S}_2$, conditioned on
$\Lambda^i\left(\mathcal{T}_1\times\mathcal{S}_1\right)$ and
$\Lambda^j\left(\mathcal{T}_2\times\mathcal{S}_2\right)$, the random
variables $N^i\left(\mathcal{T}_1\times\mathcal{S}_1\right)$ and
$N^j\left(\mathcal{T}_2\times\mathcal{S}_2\right)$ are independent
for~$i=a,b$ and~$j=a,b$.

Let $(\Omega,\mathscr{F},P)$ be the underlying probability space for
the above stochastic model. Define~$\mathscr{B}_t^i$ as the
$\sigma$-algebra generated by the space-time process~$i$
over~$[0,t)$. We define the counting process $N_t^i$ as the number
of points occurred during $[0,t)$ over the entire surface of
photodetector~$i$ regardless of their location, i.e.,
$N_t^i=N^i\bigl([0,t)\times\mathbb{R}^2\bigr)$. We say~$u_t^a$
and~$u_t^b$ are admissible controls if~$u_t^i$ is
$\mathscr{B}_t^i$-measurable and the solution
to~\eqref{CooperativeStateSpaceModel} is well defined for~$i=a,b$.

Let~$T$ be an arbitrary positive constant. We define the objective
functional
\begin{equation}\label{CooperativeCostFunctional}
J=\mathrm{E}\left[\int_0^T
\left(\alpha^a\nu_t^a\exp\bigl(-\rho\|C_tx_t^b\|^2\bigr)+
\alpha^b\nu_t^b\exp\bigl(-\rho\|C_tx_t^a\|^2\bigr)\right)dt\right]
\end{equation}
where~$\alpha^i\geqslant0,~i=a,b$. Evidently,~$J$ is a linear
combination of the expected optical energy received by the stations
during~$[0,T]$. Our goal is to obtain admissible controls~$u_t^a$
and~$u_t^b$ that maximize~$J$.
%
%
\section{The Control Problem}
We first outline some results in the estimation of the state
vector~$x_t^i$ due to Rhodes and Snyder~\cite{ART.Rhodes.Snyder.77}.
Later, we will utilize these results to approach the control
problem. Before moving forward, let us fix some notation.
Let~$(t_{k-1},t_k]$ be the interval between two successive
occurrence of the space-time process and $r_k$ be the location
of~$k^{th}$ occurring point. Assume that~$h_t\left(r,\xi_t\right)$
is continuous in~$r$ and left continuous in~$t$ and~$\xi_t$. Then
the stochastic differential equation
\begin{equation*}
d\xi_t=\int_{\mathbb{R}^2}h_t\left(r,\xi_t\right)N\left(dt\times{dr}\right)
\end{equation*}
is defined such that~$d\xi_t=0$ during~$(t_{k-1},t_k]$ and~$\xi_t$
encounters a jump of $h_{t_k}\left(r_k,\xi_{t_k}\right)$ at~$t=t_k$.

Consider the state-space model~\eqref{CooperativeStateSpaceModel}
and its associated space-time observation with the rate
process~\eqref{CooperativePoissonRate}. Assume that we are given the
increasing family of $\sigma$-algebras~$\mathscr{B}_t^i$ and~$u_t^i$
is $\mathscr{B}_t^i$-measurable. Then, regardless of the nature
of~$\mu_t^i$, the posterior
density~$p_{x_t^i}\left(x|\mathscr{B}_t^i\right)$ is Gaussian with
mean~$\hat{x}_t^i$ and covariance matrix~$\Sigma_t^i$ determined
from the stochastic differential
equations~\cite{ART.Rhodes.Snyder.77}
\begin{align}\label{}
d\hat{x}_t^i&=A_t\hat{x}_t^idt+B_tu_t^idt
+\int_{\mathbb{R}^2}M_t^i\left(r-C_t\hat{x}_t^i\right)N^i\left(dt\times{dr}\right)\label{CooperativeSDEInfinit1}\\
d\Sigma_t^i&=A_t\Sigma_t^idt+\Sigma_t^iA_t^Tdt+D_tD_t^Tdt-M_t^iC_t\Sigma_t^idN_t^i
\label{CooperativeSDEInfinit2}
\end{align}
with initial states $\hat{x}_0^i=\bar{x}_0^i$ and
$\Sigma_0^i=\bar{\Sigma}_0^i$, where $M_t^i$ is defined~as
\begin{equation*}
M_t^i=\Sigma_t^iC_t^T\left(C_t\Sigma_t^iC_t^T+R_t\right)^{-1}
\end{equation*}
Moreover, the conditional covariance matrix~$\Sigma_t^i$ obtained
from~\eqref{CooperativeSDEInfinit2} is almost surely positive
definite for~$t>0$, provided that~$\bar{\Sigma}_0^i$ is positive
definite.

Note that the formulas~\eqref{CooperativeSDEInfinit1}
and~\eqref{CooperativeSDEInfinit2} not explicitly depend
on~$\{\mu_t^i,t\geqslant0\}$, however, the estimates~$\hat{x}_t^i$
and~$\Sigma_t^i$ depend on~$\{\mu_t^i,t\geqslant0\}$ through the
observation~$N^i\left(dt\times{dr}\right)$. This dependence can be
explained by observing from~\eqref{CooperativeSDEInfinit2} that the
occurrence of each event in the space-time process subtracts the
positive definite matrix~$M_t^iC_t\Sigma_t^i$ from~$\Sigma_t^i$,
thus a larger~$\mu_t^i$ leads to a smaller estimation error
covariance by increasing the occurrence rate of the events.
According to~\eqref{CooperativeMUDefinition}, a smaller pointing
error~$\|C_tx_t^b\|$ at station~$b$ results in a larger~$\mu_t^a$
and, as a consequence, a closer estimation for~$x_t^a$, which in
turn, leads to a smaller pointing error at station~$a$. This
explains the mechanism which couples the dynamics of the stations.

We exploit the above results to prove
theorem~\ref{CooperativeMainTheorem} below which determines an upper
bound on~$J$ and establishes the conditions on~$\hat{x}_t^a$
and~$\hat{x}_t^b$ under which the upper bound can be achieved.
Before stating the theorem, we fix notation.
Let~$\Sigma=[\sigma_{ij}]$ denote a symmetric $n\times{n}$ matrix
and~$f\left(\Sigma\right)$ be a scalar function of~$\Sigma$. Assume
that the partial derivatives of $f\left(\Sigma\right)$ with respect
to elements of~$\Sigma$ exist. We denote
by~$\partial{f\left(\Sigma\right)/\partial{\Sigma}}$ a $n\times{n}$
symmetric matrix $F\left(\Sigma\right)=[F_{ij}\left(\Sigma\right)]$
such that $F_{ii}=\partial{f/\partial{\sigma_{ii}}}$ and
$F_{ij}=\left(1/2\right)\partial{f/\partial{\sigma_{ij}}}$
for~$i{\neq}j$. Let $g_t\left(\Sigma^a,\Sigma^b\right)$ be a
function of $n\times{n}$ symmetric matrices $\Sigma^a$ and
$\Sigma^b$ with values in~$\mathbb{R}$. We define the linear
operators~$\mathcal{L}_t^a\left\{\cdot\right\}$
and~$\mathcal{L}_t^b\left\{\cdot\right\}$~as
\begin{equation}\label{CooperativeLinearOperatorDefinition}
\begin{split}
\textstyle\mathcal{L}_t^a\left\{g_t\left(\Sigma^a,\Sigma^b\right)\right\}&=
g_t\textstyle\left(S_t\left(\Sigma^a\right),\Sigma^b\right)-g_t\left(\Sigma^a,\Sigma^b\right)\\
\textstyle\mathcal{L}_t^b\left\{g_t\left(\Sigma^a,\Sigma^b\right)\right\}&=
g_t\textstyle\left(\Sigma^a,S_t\left(\Sigma^b\right)\right)
-g_t\left(\Sigma^a,\Sigma^{\!\:b}\right)
\end{split}
\end{equation}
where~$S_t\left(\cdot\right)$ is defined~as
\begin{equation}\label{CooperativeDeltaDefinition}
S_t\left(\Sigma\right)=\Sigma-{\Sigma}C_t^T\bigl(C_t{\Sigma}C_t^T+R_t\bigr)^{-1}C_t\Sigma
\end{equation}
With~$I$ being the $2{\times}2$ identity matrix, we define the
$2{\times}2$ positive definite matrix~$Q_t^i$~as
\begin{equation}\label{CooperativeQMatrixDefinition}
Q_t^i=\left(I+2{\rho}C_t\Sigma_t^iC_t^T\right)^{-1/2}
\end{equation}
Also we define~$q_t^i>0$~as
\begin{equation}\label{CooperativeQDetDefinition}
q_t^i=\det\left(Q_t^i\right)
\end{equation}

\begin{theorem}\label{CooperativeMainTheorem}
Fix sample paths~$\nu_t^a$ and~$\nu_t^b$, $t\in[0,T]$.
Let~$\Sigma^a$ and~$\Sigma^b$ be $n{\times}n$ symmetric matrices.
Assume that~$g_t\left(\Sigma^a,\Sigma^b\right)$, $t\in[0,T]$ is the
backward solution to the partial differential equation
\begin{align}\label{CooperativeMainTheoremPDE}
-\frac{{\partial}g_t\left(\Sigma^a,\Sigma^b\right)}{{\partial}t}&=
\nu_t^a\frac{\alpha^a+\mathcal{L}_t^a\left\{g_t\left(\Sigma^a,\Sigma^b\right)\right\}}
{\displaystyle\sqrt{\det\left(I+2{\rho}C_t\Sigma^bC_t^T\right)}}
+\nu_t^b\frac{\alpha^b+\mathcal{L}_t^b\left\{g_t\left(\Sigma^a,\Sigma^b\right)\right\}}
{\displaystyle\sqrt{\det\left(I+2{\rho}C_t\Sigma^aC_t^T\right)}}\nonumber\\
&+\mathrm{tr}\left\{\frac{{\partial}g_t\left(\Sigma^a,\Sigma^b\right)}{{\partial}\Sigma^a}
\Bigl(A_t\Sigma^a+\Sigma^aA_t^T+D_tD_t^T\Bigr)\right\}\nonumber\\
&+\mathrm{tr}\left\{\frac{{\partial}g_t\left(\Sigma^a,\Sigma^b\right)}{{\partial}\Sigma^b}
\Bigl(A_t\Sigma^b+\Sigma^bA_t^T+D_tD_t^T\Bigr)\right\}
\end{align}
with boundary condition $g_T\left(\cdot,\cdot\right)=0$. Then for
fixed sample paths of~$\nu_t^a$ and~$\nu_t^b$, the objective
functional~\eqref{CooperativeCostFunctional} can be expressed~as
\begin{align}\label{CooperativeMainTheoremCostFunctional}
J=g_0\textstyle\left(\Sigma_0^a,\Sigma_0^b\right)
&-\mathrm{E}\left[\displaystyle\int_0^T\nu_t^aq_t^b\Bigl(\alpha^a+\mathcal{L}_t^a
\left\{g_t\textstyle\left(\Sigma^a_t,\Sigma^b_t\right)\right\}\Bigr)
\,\Bigl\{1-\exp\Bigl(-\rho\|Q_t^bC_t\hat{x}_t^b\|^2\Bigr)\Bigr\}dt\right]\nonumber\\
&-\mathrm{E}\left[\int_0^T\nu_t^bq_t^a\Bigl(\alpha^b+\mathcal{L}_t^b
\left\{g_t\textstyle\left(\Sigma^a_t,\Sigma^b_t\right)\right\}\Bigr)%
\Bigl\{1-\exp\Bigl(-\rho\|Q_t^aC_t\hat{x}_t^a\|^2\Bigr)\Bigr\}dt\right]
\end{align}
Moreover, if~$\Sigma^a$ and~$\Sigma^b$ are positive definite,
$\mathcal{L}_t^a\left\{g_t\left(\Sigma^a,\Sigma^b\right)\right\}$
and
$\mathcal{L}_t^b\left\{g_t\left(\Sigma^a,\Sigma^b\right)\right\}$
are positive for~$t\in[0,T)$, provided that~$C_t$ is full rank for
any~$t\in[0,T)$.
\end{theorem}
\begin{proof}
See appendix~\ref{CooperativeMainTheoremProof}.
\end{proof}
\begin{corollary}\label{CooperativeCorollary1}
For fixed sample paths~$\nu_t^a$ and~$\nu_t^b$ and any choice of
positive definite matrices~$\Sigma_0^a$ and~$\Sigma_0^b$,
$J^*{\triangleq\:}g_0\left(\Sigma_0^a,\Sigma_0^b\right)$ is an upper
bound for~$J$, i.e., $J\leqslant{J^*}$, and equality holds if and
only if~$C_t\hat{x}_t^a=0$ for~$t\in\mathscr{T}^a$ and
$C_t\hat{x}_t^b=0$ for~$t\in\mathscr{T}^b$ almost everywhere.
Here~$\mathscr{T}^a$ and~$\mathscr{T}^b$ are defined for fixed
sample paths~$\nu_t^a$ and~$\nu_t^b$ as
$\mathscr{T}^a=\{t\;\!|\;\!\nu_t^b\neq0,\,t\in[0,T)\}$ and
$\mathscr{T}^b=\{t\;\!|\;\!\nu_t^a\neq0,\,t\in[0,T)\}$.
\end{corollary}
\begin{proof}
The second statement of theorem~\ref{CooperativeMainTheorem} in
conjunction with positive definiteness of~$\Sigma_t^i,~i=a,b$
results
$\mathcal{L}_t^i\left\{g_t\left(\Sigma_t^a,\Sigma_t^b\right)\right\}>0$
for $i=a,b$ and $t\in[0,T)$. It follows that the integrands in the
first and second integrals
of~\eqref{CooperativeMainTheoremCostFunctional} are positive
over~$\mathscr{T}^b$ and~$\mathscr{T}^a$, respectively, except
when~$C_t\hat{x}_t^b=0$ for~$t\in\mathscr{T}^b$ and
$C_t\hat{x}_t^a=0$ for~$t\in\mathscr{T}^a$, in which they are equal
to zero. This leads to $J-J^*\leqslant0$ with equality when the
integrals vanish. The last condition holds if and only
if~$C_t\hat{x}_t^b=0$ for~$t\in\mathscr{T}^b$ and $C_t\hat{x}_t^a=0$
for~$t\in\mathscr{T}^a$ almost everywhere.
\end{proof}
\begin{corollary}\label{CooperativeCorollary2}
For nonnegative stochastic processes~$\nu_t^a$ and~$\nu_t^b$ with
piecewise continuous sample paths and nonzero expectations, the
objective functional~\eqref{CooperativeCostFunctional} achieves its
maximum if and only~if
\begin{equation}\label{CooperativeCorollary2Condition}
C_t\hat{x}_t^a=C_t\hat{x}_t^b=0
\end{equation}
for~$t\in[0,T)$ almost everywhere.
\end{corollary}
\begin{proof}
If~\eqref{CooperativeCorollary2Condition} holds,
corollary~\ref{CooperativeCorollary1} implies that for any sample
path of~$\left(\nu_t^a,\nu_t^b\right)$ the objective functional~$J$
associated to that sample path meets its upper bound. This suggests
that~\eqref{CooperativeCorollary2Condition} is a sufficient
condition for~$J$ to achieve its maximum. To show
that~\eqref{CooperativeCorollary2Condition} is a necessary
condition, assume that for some interval
$\mathscr{I}\subseteq[0,T)$, $C_t\hat{x}_t^a\neq0$ or
$C_t\hat{x}_t^b\neq0$. Because~$\nu_t^a$ and~$\nu_t^b$ have nonzero
expectations, with nonzero probability some of their sample paths
are positive over~$\mathscr{I}$. Then,
corollary~\ref{CooperativeCorollary1} implies that for those sample
paths the condition for achieving the maximum is not satisfied.
Therefore, the objective functional~$J$ cannot achieve its maximum.
\end{proof}
\begin{remark}
According to~\eqref{CooperativeSDEInfinit1}
and~\eqref{CooperativeSDEInfinit2}, $\hat{x}_t^a$ and~$\hat{x}_t^b$
do not explicitly depend on~$\nu_t^a$ and~$\nu_t^b$. This suggests
that the optimal control law not explicitly depends on~$\nu_t^a$
and~$\nu_t^b$.
\end{remark}
\begin{remark}
The condition which leads to the upper bound~$J^*$ does not depend
on~$\alpha^a$ and~$\alpha^b$. In particular, the condition is same
for $(\alpha^a,\alpha^b)=(1,0)$ and $(\alpha^a,\alpha^b)=(0,1)$.
This means that if $C_t\hat{x}_t^a=C_t\hat{x}_t^b=0$ holds during
$t\in[0,T)$, both stations receive the maximum possible optical
energy.
\end{remark}

The following lemma proposes a control law which leads to the
condition~$C_t\hat{x}_t^a=C_t\hat{x}_t^b=0$.
\begin{lemma}\label{CooperativeOptimalLemmaControl}
Consider the stochastic dynamical
system~\eqref{CooperativeSDEInfinit1} and assume
that~$C_0\hat{x}_0^i=0$. Let~$C_tB_t$ be nonsingular and~$C_t$ be
differentiable for $t\geqslant0$. Then under control
\begin{equation}\label{CooperativeOptimalControl}
u_t^idt=-\left(C_tB_t\right)^{-1}\left\{\left(C_tA_t+\dot{C}_t\right)\hat{x}_t^idt
+C_tM_t^i\int_{\mathbb{R}^2}rN^i\left(dt\times{dr}\right)\right\}
\end{equation}
we have $C_t\hat{x}_t^i=0$ for any $t\geqslant0$.
\end{lemma}
\begin{proof}
We verify the validity of the lemma by putting
\begin{equation}\label{CooperativeOptimalControlLemmaProof}
u_t^idt=-\left(C_tB_t\right)^{-1}\left\{\left(C_tA_t+\dot{C}_t\right)\hat{x}_t^idt
+\int_{\mathbb{R}^2}C_tM_t^i\left(r-C_t\hat{x}_t^i\right)N^i\left(dt\times{dr}\right)\right\}
\end{equation}
into~\eqref{CooperativeSDEInfinit1} and left multiplying both sides
by~$C_t$. The resulting equation will be
$C_td\hat{x}_t^i=-\dot{C}_t\hat{x}_t^idt$, which yields
$d\left(C_t\hat{x}_t^i\right)=0$. Then we argue that
$C_t\hat{x}_t^i=C_0\hat{x}_0^i=0$ for~$t\geqslant0$. Finally, we
put~$C_t\hat{x}_t^i=0$
into~\eqref{CooperativeOptimalControlLemmaProof} to
obtain~\eqref{CooperativeOptimalControl}.
\end{proof}
%
%
\section{Conclusion}
Cooperative optical beam tracking, a scheme for maintaining
alignment in a free-space optical link, has been addressed. A
stochastic model has been developed which captures three sources of
randomness: relative motion of stations, characteristic of
photodetectors, and fluctuation of optical power caused by optical
fade and information-bearing signals modulating the optical beams.
An optimal control law has been proposed which maximizes the
expected optical energy received by stations of the link. It has
been shown that under moderate assumptions, the control law does not
depend on the characteristic of the optical fade or
information-bearing signals modulating the optical beams.
%
%
\appendix
\section{Proof of Theorem~\ref{CooperativeMainTheorem}}\label{CooperativeMainTheoremProof}
We prove the first statement of the theorem through the following
four steps.

\noindent Step I: Recalling that
$p_{x_t^i}\bigl(x|\mathscr{B}_t^i\bigr)$ is Gaussian with
mean~$\hat{x}_t^i$ and covariance~$\Sigma_t^i$, it is
straightforward to show that
\begin{align*}
\mathrm{E}\left[\exp\bigl(-\rho\|C_tx_t^i\|^2\bigr)\right]&=
\mathrm{E}\left[\mathrm{E}\left[\exp\bigl(-\rho\|C_tx_t^i\|^2\bigr)\big{|}\mathscr{B}_t^i\right]\right]\\
&=\mathrm{E}\left[q_t^i\exp\bigl(-\rho\|Q_t^iC_t\hat{x}_t^i\|^2\bigr)\right]
\end{align*}
where~$Q_t^i$ and~$q_t^i$ are defined
by~\eqref{CooperativeQMatrixDefinition}
and~\eqref{CooperativeQDetDefinition}, respectively. Using the above
equation we rewrite~\eqref{CooperativeCostFunctional}~as
\begin{align}\label{CooperativeMainTheoremProofEQ2}
J=&\mathrm{E}\left[\displaystyle\int_0^T\textstyle\left(\alpha^a\nu_t^aq_t^b
+\alpha^b\nu_t^bq_t^a\right)dt\right]\nonumber\\
-&\mathrm{E}\left[\displaystyle\int_0^T\alpha^a\nu_t^aq_t^b
\left\{1-\exp\textstyle\left(-\rho\|Q_t^bC_t\hat{x}_t^b\|^2\right)\right\}dt\right]\nonumber\\
-&\mathrm{E}\left[\displaystyle\int_0^T\alpha^b\nu_t^bq_t^a
\left\{1-\exp\textstyle\left(-\rho\|Q_t^aC_t\hat{x}_t^a\|^2\right)\right\}dt\right]
\end{align}

\noindent Step II: Let $(i,j)=(a,b),(b,a)$. For any $t\geqslant0$
and for any small positive~$\epsilon$,
${\Delta}N_t^i{\triangleq}N_{t+\epsilon}^i-N_t^i$ is a Poisson
random variable conditioned on the rate
\begin{equation*}
\Lambda_t^i=\int_t^{t+\epsilon}\nu_\tau^i\exp\bigl(-\rho\|C_{\tau}x_\tau^j\|^2\bigr)d\tau
\end{equation*}
Thus, using the law of total probability we can write
\begin{align}\label{CooperativeMainTheoremProofEQ3}
\Pr\left\{{\Delta}N_t^i=1\big{|}\mathscr{B}_t^j\right\}
&=\mathrm{E}\left[\Pr\left\{{\Delta}N_t^i=1|\Lambda_t^i\right\}\big{|}\mathscr{B}_t^j\right]\nonumber\\
&=\mathrm{E}\left[\Lambda_t^ie^{-\Lambda_t^i}\big{|}\mathscr{B}_t^j\right]\nonumber\\
&={p_{t,\epsilon}^i}+O\left(\epsilon^2\right)
\end{align}
where~$p_{t,\epsilon}^i$ is defined~as
\begin{equation}\label{CooperativeMainTheoremProofEQ4}
p_{t,\epsilon}^i=\int_t^{t+\epsilon}\nu_\tau^iq_\tau^j\exp\bigl(-\rho\|Q_\tau^jC_{\tau}x_\tau^j\|^2\bigr)d\tau
\end{equation}
In a similar manner, we can show that
\begin{equation}\label{CooperativeMainTheoremProofEQ5}
\begin{split}
\Pr\left\{{\Delta}N_t^i=0\big{|}\mathscr{B}_t^j\right\}&=1-{p_{t,\epsilon}^i}+O\left(\epsilon^2\right)\\
\Pr\left\{{\Delta}N_t^i\geqslant2\big{|}\mathscr{B}_t^j\right\}&=O\left(\epsilon^2\right)
\end{split}
\end{equation}

\noindent Step III: Let $f\left(\Sigma\right)$ be a scalar function
of $n{\times}n$ symmetric matrix~$\Sigma$ and assume~$f$ is
differentiable with respect to~$\Sigma$. We can write
\begin{align}\label{CooperativeMainTheoremProofEQ6}
\mathrm{E}\left[f\left(\Sigma_{t+\epsilon}^i\right)\right]&=
\mathrm{E}\left[\mathrm{E}\left[f\left(\Sigma_{t+\epsilon}^i\right)\big{|}\mathscr{B}_t^j\right]\right]\nonumber\\
&=\mathrm{E}\left[\:\sum_{k=0}^\infty\mathrm{E}\left[f\left(\Sigma_{t+\epsilon}^i\right)
\big{|}{\Delta}N_t^i=k,\mathscr{B}_t^j\right]
\Pr\left\{{\Delta}N_t^i=k\big{|}\mathscr{B}_t^j\right\}\right]\nonumber\\
&=\mathrm{E}\left[\:\sum_{k=0}^1\mathrm{E}\left[f\left(\Sigma_{t+\epsilon}^i\right)
\big{|}{\Delta}N_t^i=k,\mathscr{B}_t^j\right]
\Pr\left\{{\Delta}N_t^i=k\big{|}\mathscr{B}_t^j\right\}\right]+O\left(\epsilon^2\right)
\end{align}
Also we have
\begin{align*}
\mathrm{E}\left[f\left(\Sigma_{t+\epsilon}^i\right)
\big{|}{\Delta}N_t^i=k,\mathscr{B}_t^j\right]&
=\mathrm{E}\biggl[\mathrm{E}\left[f\left(\Sigma_{t+\epsilon}^i\right)
\big{|}\Sigma_t^i,{\Delta}N_t^i=k,\mathscr{B}_t^j\right]
\Big{|}k,\mathscr{B}_t^j\biggr]\\
&=\mathrm{E}\biggl[\mathrm{E}\left[f\left(\Sigma_{t+\epsilon}^i\right)
\big{|}\Sigma_t^i,{\Delta}N_t^i=k\right]\Big{|}k,\mathscr{B}_t^j\biggr]
\end{align*}
Inserting the last equation
into~\eqref{CooperativeMainTheoremProofEQ6}, after some manipulation
we find
\begin{equation}\label{CooperativeMainTheoremProofEQ7}
\mathrm{E}\left[f\left(\Sigma_{t+\epsilon}^i\right)\right]=
\mathrm{E}\left[\:\sum_{k=0}^1\mathrm{E}\left[f\left(\Sigma_{t+\epsilon}^i\right)
\big{|}\Sigma_t^i,{\Delta}N_t^i=k\right]\Pr\left\{{\Delta}N_t^i=k
\big{|}\mathscr{B}_t^j\right\}\right]+O\left(\epsilon^2\right)
\end{equation}
Conditioned on~$\Sigma_t^i$ and ${\Delta}N_t^i=0,1$, the stochastic
differential equation~\eqref{CooperativeSDEInfinit2} can be solved
for~$\Sigma_{t+\epsilon}^i$. This solution leads~to
\begin{align*}
\mathrm{E}\left[f\left(\Sigma_{t+\epsilon}^i\right)\big{|}\Sigma_t^i,{\Delta}N_t^i=0\right]
&=f\left(\Sigma_t^i+{\epsilon}A_t\Sigma_t^i+{\epsilon}\Sigma_t^iA_t^T
+{\epsilon}D_tD_t^T+O\left(\epsilon^2\right)\right)\\
\mathrm{E}\left[f\left(\Sigma_{t+\epsilon}^i\right)\big{|}\Sigma_t^i,{\Delta}N_t^i=1\right]&=
f\left(\Sigma_t^i-M_t^iC_t\Sigma_t^i+O\left(\epsilon\right)\right)
\end{align*}
We linearize the above equations with respect to~$\epsilon$ and plug
the linearized equations together
with~\eqref{CooperativeMainTheoremProofEQ3}
and~\eqref{CooperativeMainTheoremProofEQ5}
into~\eqref{CooperativeMainTheoremProofEQ7} to obtain
\begin{align}\label{CooperativeMainTheoremProofEQ8}
\mathrm{E}\left[f\left(\Sigma_{t+\epsilon}^i\right)\right]&=
\mathrm{E}\Big[f\left(\Sigma_t^i\right)+{\epsilon}\,\mathrm{tr}
\Bigl\{\left({\partial}f\left(\Sigma_t^i\right)/{\partial}\Sigma_t^i\right)
\left(A_t\Sigma_t^i+\Sigma_t^iA_t^T+D_tD_t^T\right)\Bigr\}\nonumber\\
&+p_{t,\epsilon}^i\Bigl(f\left(\Sigma_t^i-M_t^iC_t\Sigma_t^i\right)
-f\left(\Sigma_t^i\right)\Bigr)\Big]+O\left(\epsilon^2\right)
\end{align}

Let $g_t\left(\Sigma^a,\Sigma^b\right)$ be a scalar function of
$n{\times}n$ symmetric matrices~$\Sigma^a$ and~$\Sigma^b$ and assume
its partial derivatives with respect to~$t$,~$\Sigma^a$
and~$\Sigma^b$ exist. Then we can write
\begin{equation*}
\mathrm{E}\left[g_{t+\epsilon}\bigl(\Sigma_{t+\epsilon}^a,\Sigma_{t+\epsilon}^b\bigr)\right]=
\mathrm{E}\left[g_t\bigl(\Sigma_{t+\epsilon}^a,\Sigma_{t+\epsilon}^b\bigr)
+{\epsilon}{\partial}g_t\bigl(\Sigma_{t+\epsilon}^a,\Sigma_{t+\epsilon}^b\bigr)/{\partial}t\right]
+O\left(\epsilon^2\right)
\end{equation*}
Applying~\eqref{CooperativeMainTheoremProofEQ8} to~$g_t$
and~$\partial{g_t}/\partial{t}$ first for~$i=a$ and then for~$i=b$,
we find
\begin{align}\label{CooperativeMainTheoremProofEQ8P}
\mathrm{E}\left[g_{t+\epsilon}\bigl(\Sigma_{t+\epsilon}^a,\Sigma_{t+\epsilon}^b\bigr)\right]
&=\mathrm{E}\Bigg[g_t\bigl(\Sigma_t^a,\Sigma_t^b\bigr)+{\epsilon}
{\partial}g_t\bigl(\Sigma_t^a,\Sigma_t^b\bigr)/{\partial}t\nonumber\\
&+{\epsilon}\,\mathrm{tr}\left\{\left({\partial}g_t\bigl(\Sigma_t^a,\Sigma_t^b\bigr)/\partial\Sigma_t^a\right)
\Bigl(A_t\Sigma_t^a+\Sigma_t^aA_t^T+D_tD_t^T\Bigr)\right\}\nonumber\\
&+{\epsilon}\,\mathrm{tr}\left\{\left({\partial}g_t\bigl(\Sigma_t^a,\Sigma_t^b\bigr)/\partial\Sigma_t^b\right)
\Bigl(A_t\Sigma_t^b+\Sigma_t^bA_t^T+D_tD_t^T\Bigr)\right\}\nonumber\\
&+p_{t,\epsilon}^a\left(g_t\bigl(\Sigma_t^a-M_t^aC_t\Sigma_t^a,\Sigma_t^b\bigr)
-g_t\bigl(\Sigma_t^a,\Sigma_t^b\bigr)\right)\nonumber\\
&+p_{t,\epsilon}^b\left(g_t\bigl(\Sigma_t^a,\Sigma_t^b-M_t^bC_t\Sigma_t^b\bigr)
-g_t\bigl(\Sigma_t^a,\Sigma_t^b\bigr)\right)\Bigg]+O\left(\epsilon^2\right)
\end{align}
In this equation, the term involving~$p_{t,\epsilon}^a$ can be
simplified~as
\begin{align*}
&\mathrm{E}\left[p_{t,\epsilon}^a\left(g_t\bigl(\Sigma_t^a-M_t^aC_t\Sigma_t^a,\Sigma_t^b\bigr)
-g_t\bigl(\Sigma_t^a,\Sigma_t^b\bigr)\right)\right]\\
&\qquad
=\mathrm{E}\left[\int_t^{t+\epsilon}{\!\!}\nu_\tau^aq_\tau^b\exp\bigl(-\rho\|Q_\tau^bC_{\tau}x_\tau^b\|^2\bigr)d\tau
\cdot\mathcal{L}_t^a\left\{g_t\textstyle\left(\Sigma_t^a,\Sigma_t^b\right)\right\}\right]\\
&\qquad=\mathrm{E}\left[\epsilon{\:\!}\nu_t^aq_t^b\exp\bigl(-\rho\|Q_t^bC_tx_t^b\|^2\bigr)
\mathcal{L}_t^a\left\{g_t\textstyle\left(\Sigma_t^a,\Sigma_t^b\right)\right\}\right]+O\left(\epsilon^2\right)
\end{align*}
Here, the first equation is obtained
from~\eqref{CooperativeLinearOperatorDefinition}
and~\eqref{CooperativeMainTheoremProofEQ4} and the second one is
concluded from the fact that with probability
$1-O\left(\epsilon\right)$ the integrand in the first equation is
continuous. Applying this result
to~\eqref{CooperativeMainTheoremProofEQ8P} and rearranging the
equation we obtain
\begin{align}\label{CooperativeMainTheoremProofEQ9}
&\mathrm{E}\left[g_{t+\epsilon}\bigl(\Sigma_{t+\epsilon}^a,\Sigma_{t+\epsilon}^b\bigr)
+\epsilon\left(\alpha^a\nu_t^aq_t^b+\alpha^b\nu_t^bq_t^a\right)\right]
-\mathrm{E}\left[g_t\bigl(\Sigma_t^a,\Sigma_t^b\bigr)-\epsilon\Gamma_t\right]+O\left(\epsilon^2\right)\nonumber\\
&\qquad=\epsilon\mathrm{E}\Bigg[\frac{{\partial}g_t\left(\Sigma_t^a,\Sigma_t^b\right)}{{\partial}t}+
\nu_t^a\frac{\alpha^a+\mathcal{L}_t^a\left\{g_t\left(\Sigma_t^a,\Sigma_t^b\right)\right\}}
{\displaystyle\sqrt{\det\left(I+2{\rho}C_t\Sigma_t^bC_t^T\right)}}
+\nu_t^b\frac{\alpha^b+\mathcal{L}_t^b\left\{g_t\left(\Sigma_t^a,\Sigma_t^b\right)\right\}}
{\displaystyle\sqrt{\det\left(I+2{\rho}C_t\Sigma_t^aC_t^T\right)}}\nonumber\\
&\qquad+\mathrm{tr}\left\{\frac{{\partial}g_t\left(\Sigma_t^a,\Sigma_t^b\right)}{{\partial}\Sigma_t^a}
\Bigl(A_t\Sigma_t^a+\Sigma_t^aA_t^T+D_tD_t^T\Bigr)\right\}\nonumber\\
&\qquad+\mathrm{tr}\left\{\frac{{\partial}g_t\left(\Sigma_t^a,\Sigma_t^b\right)}{{\partial}\Sigma_t^b}
\Bigl(A_t\Sigma_t^b+\Sigma_t^bA_t^T+D_tD_t^T\Bigr)\right\}\Bigg]
\end{align}
where~$\Gamma_t$ is defined~as
\begin{align}\label{CooperativeMainTheoremProofEQ10}
\Gamma_t&=\nu_t^aq_t^b\mathcal{L}_t^a
\left\{g_t\textstyle\left(\Sigma^a_t,\Sigma^b_t\right)\right\}
\left\{1-\exp\textstyle\left(-\rho\|Q_t^bC_t\hat{x}_t^b\|^2\right)\right\}\nonumber\\
&+\nu_t^bq_t^a\mathcal{L}_t^b
\left\{g_t\textstyle\left(\Sigma^a_t,\Sigma^b_t\right)\right\}%
\left\{1-\exp\textstyle\left(-\rho\|Q_t^aC_t\hat{x}_t^a\|^2\right)\right\}
\end{align}
Assuming that~$g_t\left(\cdot,\cdot\right)$ is the backward solution
to~\eqref{CooperativeMainTheoremPDE} with boundary
condition~$g_T\left(\cdot,\cdot\right)=0$, the right hand side
of~\eqref{CooperativeMainTheoremProofEQ9} is identically zero, which
yields~to
\begin{equation}\label{CooperativeMainTheoremProofEQ11}
\mathrm{E}\left[g_{t+\epsilon}\bigl(\Sigma_{t+\epsilon}^a,\Sigma_{t+\epsilon}^b\bigr)+\epsilon
\textstyle\left(\alpha^a\nu_t^aq_t^b+\alpha^b\nu_t^bq_t^a\right)\right]
=\mathrm{E}\left[g_t\bigl(\Sigma_t^a,\Sigma_t^b\bigr)-\epsilon\Gamma_t\right]+O\left(\epsilon^2\right)
\end{equation}

\noindent Step IV: Let us denote the first term in the right hand
side of~\eqref{CooperativeMainTheoremProofEQ2} by~$\tilde{J}$. We
partition the interval $[0,T)$ into $K$ subintervals
$[t_k,t_{k+1}),~k=0,1,\ldots,K-1$, where $t_0=0$, $t_K=T$, and
$t_{k+1}-t_k\triangleq\epsilon_k>0$. Then,~$\tilde{J}$ can be
approximated~by
\begin{equation*}
\tilde{J}{\simeq}\tilde{J}_K=\sum_{k=0}^{K-1}\epsilon_k\mathrm{E}\left[\alpha^a\nu_{t_k}^aq_{t_k}^b
+\alpha^b\nu_{t_k}^bq_{t_k}^a\right]
+\mathrm{E}\left[g_{t_K}\bigl(\Sigma_{t_K}^a,\Sigma_{t_K}^b\bigr)\right]
\end{equation*}
Note that
$g_{t_K}\left(\cdot,\cdot\right)=g_T\left(\cdot,\cdot\right)=0$,
thus the last term does not affect the sum and is arbitrarily added
to the right hand side. The above equation can be rewritten~as
\begin{align*}
\tilde{J}_K&=\sum_{k=0}^{K-2}\epsilon_k\mathrm{E}\left[\alpha^a\nu_{t_k}^aq_{t_k}^b
+\alpha^b\nu_{t_k}^bq_{t_k}^a\right]\\
&+\mathrm{E}\left[g_{t_K}\bigl(\Sigma_{t_K}^a,\Sigma_{t_K}^b\bigr)
+\epsilon_{K-1}\left(\alpha^a\nu_{t_{K-1}}^aq_{t_{K-1}}^b
+\alpha^b\nu_{t_{K-1}}^bq_{t_{K-1}}^a\right)\right]
\end{align*}
Since $t_K=t_{K-1}+\epsilon_{K-1}$, we can
apply~\eqref{CooperativeMainTheoremProofEQ11} to the second term of
the right hand side of this equation to obtain
\begin{align*}
\tilde{J}_K&=\sum_{k=0}^{K-2}\epsilon_k\mathrm{E}\left[\alpha^a\nu_{t_k}^aq_{t_k}^b
+\alpha^b\nu_{t_k}^bq_{t_k}^a\right]\\
&+\mathrm{E}\left[g_{t_{K-1}}\bigl(\Sigma_{t_{K-1}}^a,\Sigma_{t_{K-1}}^b\bigr)\right]
-\epsilon_{K-1}\mathrm{E}\left[\Gamma_{t_{K-1}}\right]+O\left(\epsilon^2_{K-1}\right)
\end{align*}
Continuing this procedure for $k=K-2,K-3,\ldots,0$, we obtain
\begin{equation*}
\tilde{J}_K=\mathrm{E}\left[g_{t_0}\left(\Sigma_{t_0}^a,\Sigma_{t_0}^b\right)\right]
-\sum_{k=0}^{K-1}\epsilon_k\mathrm{E}\left[\Gamma_{t_k}\right]
+\sum_{k=0}^{K-1}O\left(\epsilon_k^2\right)
\end{equation*}
Now we let $K\to\infty$ and $\displaystyle\max\epsilon_k\to0$ to
obtain
\begin{equation*}
\tilde{J}=\lim_{\substack{N\to\infty\\\max\epsilon_k\to0}}\tilde{J}_K
=g_0\bigl(\Sigma_0^a,\Sigma_0^b\bigr)-\mathrm{E}\left[\int_0^T\Gamma_tdt\right]
\end{equation*}
Finally, we put this equation
into~\eqref{CooperativeMainTheoremProofEQ2} to
obtain~\eqref{CooperativeMainTheoremCostFunctional}.

In order to prove the second statement of the theorem, we need the
following preliminaries:
\begin{enumerate}
\item [P-1)] In the context of this proof, we say
$f\left(\cdot\right):\mathbb{R}^{n\times{n}}\to\mathbb{R}$ is
strictly decreasing if for any symmetric positive definite
matrices~$\Sigma$ and~$\Delta$, we have
$f\left(\Sigma+\Delta\right)<f\left(\Sigma\right)$. Also we say
$f\left(\cdot\right)$ is m-positive if for any symmetric positive
definite matrix~$\Sigma$ we have $f\left(\Sigma\right)>0$.
\item [P-2)]
If~$f_1\left(\cdot\right)$ and~$f_2\left(\cdot\right)$ are strictly
decreasing and m-positive,
$f_1\left(\cdot\right)+f_2\left(\cdot\right)$ and
$f_1\left(\cdot\right)f_2\left(\cdot\right)$ are strictly decreasing
and m-positive as well.
\item [P-3)] If~$f\left(\cdot\right)$ is
strictly decreasing, for any~$t$ in which~$C_t$ is full rank and any
positive definite~$\Sigma$ we have
$f\left(\Sigma\right)<f\left(S_t\left(\Sigma\right)\right)$,
where~$S_t\left(\cdot\right)$ is defined
by~\eqref{CooperativeDeltaDefinition}.
\begin{proof}
Applying the matrix inversion lemma
to~\eqref{CooperativeDeltaDefinition}, it is easy to verify that
$S_t\left(\Sigma\right)$ is positive definite for any positive
definite~$\Sigma$. Also we know
from~\eqref{CooperativeDeltaDefinition} that when~$C_t$ is full
rank, $\Delta\triangleq\Sigma-S_t\left(\Sigma\right)$ is a positive
definite matrix. Because $f\left(\cdot\right)$ is assumed to be
strictly decreasing, we can write
$f\left(\Sigma\right)=f\left(S_t\left(\Sigma\right)+\Delta\right)<f\left(S_t\left(\Sigma\right)\right)$.
\end{proof}
\item [P-4)] If~$f\left(\cdot\right)$ is strictly decreasing and
m-positive, for any fixed~$t$, $f\left(S_t\left(\cdot\right)\right)$
is strictly decreasing and m-positive.
\begin{proof}
It is easy to verify that
\begin{equation}\label{CooperativeMainTheoremProofEQ20}
\Sigma^{-1}-\left(\Sigma+\Delta\right)^{-1}=\left(\Sigma+\Sigma\Delta^{-1}\Sigma\right)^{-1}
\end{equation}
holds for any invertible matrices~$\Sigma$ and~$\Delta$.
Let~$\Sigma$ and~$\Delta$ be positive definite matrices.
Then~\eqref{CooperativeMainTheoremProofEQ20} implies that
$\tilde{\Delta}{\,\triangleq\,}\Sigma^{-1}-\left(\Sigma+\Delta\right)^{-1}$
is a positive definite matrix. Using the matrix inversion lemma and
replacing $\Sigma^{-1}$ with
$\left(\Sigma+\Delta\right)^{-1}+\tilde{\Delta}$, we can write
\begin{align*}
&S_t\left(\Sigma+\Delta\right)-S_t\left(\Sigma\right)\\
&\qquad\quad=\left(\left(\Sigma+\Delta\right)^{-1}+C_t^TR_t^{-1}C_t\right)^{-1}
-\left(\Sigma^{-1}+C_t^TR_t^{-1}C_t\right)^{-1}\\
&\qquad\quad=\left(\left(\Sigma+\Delta\right)^{-1}+C_t^TR_t^{-1}C_t\right)^{-1}
-\left(\left(\left(\Sigma+\Delta\right)^{-1}+C_t^TR_t^{-1}C_t\right)+\tilde{\Delta}\right)^{-1}
\end{align*}
Applying identity~\eqref{CooperativeMainTheoremProofEQ20} to the
last equation, we find that
$S_t\left(\Sigma+\Delta\right)-S_t\left(\Sigma\right)$ is positive
definite. Then, because $S_t\left(\Sigma\right)$ is positive
definite and~$f\left(\cdot\right)$ is strictly decreasing, we have
\begin{equation*}
f\left(S_t\left(\Sigma+\Delta\right)\right)=f\left(S_t\left(\Sigma\right)
+\left\{S_t\left(\Sigma+\Delta\right)-S_t\left(\Sigma\right)\right\}\right)<f\left(S_t\left(\Sigma\right)\right)
\end{equation*}
which means $f\left(S_t\left(\cdot\right)\right)$ is strictly
decreasing. Moreover, because~$f\left(\cdot\right)$ is m-positive
and~$S_t\left(\Sigma\right)$ is positive definite,
$f\left(S_t\left(\cdot\right)\right)$ is m-positive.
\end{proof}
\item [P-5)] For any fixed~$t$ in which~$C_t$ is full rank,~$h_t\left(\Sigma\right)$ defined~as
\begin{equation*}
h_t\left(\Sigma\right)=\frac{1}{\sqrt{\det\left(I+2{\rho}C_t{\Sigma}C_t^T\right)}}
\end{equation*}
is strictly decreasing and m-positive.
\begin{proof}
For positive definite~$\Sigma$ and~$\Delta$ we can write
\begin{align*}
\frac{h_t\left(\Sigma\right)}{h_t\left(\Sigma+\Delta\right)}&=
\sqrt{\frac{\det\left(I+2{\rho}C_t{\Sigma}C_t^T+2{\rho}C_t{\Delta}C_t^T\right)}
{\det\left(I+2{\rho}C_t{\Sigma}C_t^T\right)}}\\
&=\sqrt{\det\left(I+2{\rho}\Delta^*\right)}
\end{align*}
where~$\Delta^*$ is defined~as
\begin{equation*}
\Delta^*=\left(\left(I+2{\rho}C_t{\Sigma}C_t^T\right)^{-1/2}C_t\right)
\Delta\left(\left(I+2{\rho}C_t{\Sigma}C_t^T\right)^{-1/2}C_t\right)^T
\end{equation*}
Because~$\Delta$ is positive definite and $C_t$ is full
rank,~$\Delta^*$ is positive definite. This implies that
$\det\left(I+2{\rho}\Delta^*\right)>1$, which leads to
$h_t\left(\Sigma+\Delta\right)<h_t\left(\Sigma\right)$.
\end{proof}
\item [P-6)] Let~$f_t\bigl(\Sigma^a,\Sigma^b\bigr)$ be a scalar
function of $n{\times}n$ matrices~$\Sigma^a$ and~$\Sigma^b$. Assume
that the function is strictly decreasing and m-positive in
both~$\Sigma^a$ and~$\Sigma^b$. For~$\epsilon>0$ define the linear
operator~$\mathcal{K}_t^\epsilon$~as
\begin{equation*}
\mathcal{K}_t^{\epsilon}f_t\bigl(\Sigma^a,\Sigma^b\bigr)=
\left(1-\epsilon\nu_t^ah_t\bigl(\Sigma^b\bigr)-\epsilon\nu_t^bh_t\bigl(\Sigma^a\bigr)\right)
f_t\bigl(X_t^\epsilon\bigl(\Sigma^a\bigr),X_t^\epsilon\bigl(\Sigma^b\bigr)\bigr)
\end{equation*}
where
\begin{equation*}
X_t^\epsilon\left(\Sigma\right)=\Sigma+\epsilon\left(A_t\Sigma+{\Sigma}A_t^T+D_tD_t^T\right)
\end{equation*}
Then, for any symmetric positive definite
matrices~$\Sigma^a$,~$\Sigma^b$,~$\Delta^a$, and~$\Delta^b$, there
exists
$\zeta=\zeta\left(\Sigma^a,\Sigma^b,\Delta^a,\Delta^b\right)>0$ such
that for any $0\leqslant\epsilon<\zeta$ we have
\begin{align*}
\mathcal{K}_t^{\epsilon}f_t\bigl(\Sigma^a+\Delta^a,\Sigma^b\bigr)&<
\mathcal{K}_t^{\epsilon}f_t\bigl(\Sigma^a,\Sigma^b\bigr)\\
\mathcal{K}_t^{\epsilon}f_t\bigl(\Sigma^a,\Sigma^b+\Delta^b\bigr)&<
\mathcal{K}_t^{\epsilon}f_t\bigl(\Sigma^a,\Sigma^b\bigr)\\
\mathcal{K}_t^{\epsilon}f_t\bigl(\Sigma^a,\Sigma^b\bigr)&>0
\end{align*}
Therefore, as $\epsilon\to0^+$, the above condition is satisfied for
any choice of $\Sigma^a$, $\Sigma^b$, $\Delta^a$, and $\Delta^b$.
This means that if~$\epsilon$ lays in a neighborhood of~$0$, however
small, $\mathcal{K}_t^{\epsilon}f_t\bigl(\Sigma^a,\Sigma^b\bigr)$ is
strictly decreasing and m-positive in both~$\Sigma^a$
and~$\Sigma^b$.
\end{enumerate}
We claim that $g_t\bigl(\Sigma^a,\Sigma^b\bigr)$, the solution
of~\eqref{CooperativeMainTheoremPDE} with boundary condition
$g_T\left(\cdot,\cdot\right)=0$, is strictly decreasing in
both~$\Sigma^a$ and~$\Sigma^b$ for any $t\in[0,T)$, provided
that~$C_t$ is full rank for $t\in[0,T)$. Once the claim is proven,
we apply (P-3) to~\eqref{CooperativeLinearOperatorDefinition} in
order to show
$\mathcal{L}_t^i\left\{g_t\left(\Sigma^a,\Sigma^b\right)\right\}>0,~i=a,b$
for any positive definite matrices~$\Sigma^a$ and~$\Sigma^b$ and
$t\in[0,T)$.

To prove our claim, for any~$0\leqslant{t}<T$, we partition the
interval $[t,T)$ into~$K$ subintervals
$[t_{k+1},t_k),~k=0,1,\ldots,K-1$, where $t_K=t$, $t_0=T$, and
$t_k-t_{k+1}=\epsilon_k>0$. It is straightforward to discretise the
partial differential equation~\eqref{CooperativeMainTheoremPDE} over
this partition to obtain the recursive equation
\begin{align}\label{CooperativeMainTheoremProofEQ12}
g_{t_{k+1}}\bigl(\Sigma^a,\Sigma^b\bigr)&=\epsilon_k\left(\alpha^a\nu_{t_k}^ah_{t_k}\bigl(\Sigma^b\bigr)
+\alpha^b\nu_{t_k}^bh_{t_k}\bigl(\Sigma^a\bigr)\right)\nonumber\\
&+\epsilon_k\left(\nu_{t_k}^ah_{t_k}\bigl(\Sigma^b\bigr)g_{t_k}\bigl(S_{t_k}\left(\Sigma^a\right),\Sigma^b\bigr)
+\nu_{t_k}^bh_{t_k}\bigl(\Sigma^a\bigr)g_{t_k}
\bigl(\Sigma^a,S_{t_k}\textstyle\left(\Sigma^b\right)\bigr)\right)\nonumber\\
&+\mathcal{K}_{t_k}^{\epsilon_k}g_{t_k}\bigl(\Sigma^a,\Sigma^b\bigr)+O\left(\epsilon_k^2\right)
\end{align}
Starting from $g_{t_0}\left(\cdot,\cdot\right)=0$ and using this
recursive equation for $k=0,1,2,\ldots{K-1}$, we can determine
$g_{t_K}\left(\cdot,\cdot\right)$. Then by letting $K\to\infty$ such
that $\max\epsilon_k\to0$, we have
$g_{t_K}\left(\cdot,\cdot\right){\to}g_t\left(\cdot,\cdot\right)$.

We prove by induction that as $K\to\infty$ and $\max\epsilon_k\to0$,
for $k=1,2,\ldots,K$, $g_{t_k}\left(\cdot,\cdot\right)$ is strictly
decreasing and m-positive in both~$\Sigma^a$ and~$\Sigma^b$. From
(P-2) and (P-5) we find out that
\begin{equation*}
g_{t_1}\bigl(\Sigma^a,\Sigma^b\bigr)=\epsilon_0\left(\alpha^a\nu_{t_0}^ah_{t_0}\bigl(\Sigma^b\bigr)
+\alpha^b\nu_{t_0}^bh_{t_0}\bigl(\Sigma^a\bigr)\right)
\end{equation*}
is strictly decreasing and m-positive. Now we show that if
$g_{t_k}\bigl(\Sigma^a,\Sigma^b\bigr)$ is strictly decreasing and
m-positive, $g_{t_{k+1}}\bigl(\Sigma^a,\Sigma^b\bigr)$ is strictly
decreasing and m-positive as well. For this purpose, we use (P-2,
P-5) and (P-2, P-4, P-5) respectively to show that the first and
second terms on the right hand side
of~\eqref{CooperativeMainTheoremProofEQ12} are strictly decreasing
and m-positive. Also as $\epsilon_k\to0^+$, (P-6) implies that the
third term on the right hand side
of~\eqref{CooperativeMainTheoremProofEQ12} is strictly decreasing
and m-positive. Because all three terms on the right hand side
of~\eqref{CooperativeMainTheoremProofEQ12} are strictly decreasing
and m-positive, we conclude from (P-2) that
$g_{t_{k+1}}\bigl(\Sigma^a,\Sigma^b\bigr)$ is strictly decreasing
and m-positive.
%
%

%
%
\end{document}